\documentclass[12pt,reqno]{article}

\usepackage[usenames]{color}
\usepackage{amssymb}
\usepackage{graphicx}
\usepackage{amscd}

\usepackage{color}
\usepackage{fullpage}
\usepackage{float}

\usepackage{graphics,amsmath,amssymb}
\usepackage{amsthm}
\usepackage{amsfonts}
\usepackage{latexsym}
\usepackage{epsf}
\usepackage{changepage}

\begin{document}

\hspace{1mm}

\begin{center}
\vskip 2mm{\LARGE\bf 
Tables of \\
Record Gaps Between Prime Constellations \\
}

\vskip 7mm
\large
Alexei Kourbatov\\
JavaScripter.net/math \\
{\tt akourbatov@gmail.com}
\end{center}

\vskip 6mm

\begin{abstract}
\noindent
We present tables of record gaps between densest prime constellations, or $k$-tuplets.
The tables contain all maximal gaps between prime $k$-tuplets up to $10^{15}$, for each $k\le7$.
\end{abstract}

\vskip 6mm
\noindent
{\bf Introduction.} Mathematicians have long known that there are {\em infinitely many primes}.
We also know that the average gaps between primes near $p$ are about $\log p$ 
(this follows from the {\em prime number theorem}).
However, this knowledge has not helped prove or disprove {\em Cram\'er's conjecture} \cite{cram}
that a record prime gap near $p$ is $O(\log^2{\negthinspace}p)$. 
Computations strongly suggest that $O(\log^2{\negthinspace}p)$ is
a realistic upper bound for the size of prime gaps; moreover, it even seems possible that
the prime gap ending at $p$ is {\bf\em always} less than $\log^2{\negthinspace}p$ \cite[A005250]{nicely,toes,shanks,oeis}. 
Still, the $O(\log^2{\negthinspace}p)$ upper bound for prime gaps remains just a plausible conjecture.

For {\em prime $k$-tuples}, the situation is somewhat similar~--- yet even more uncertain. 
Zhang's theorem \cite{zhang} implies that {\em at least some types of prime $k$-tuples} do occur infinitely often.
The {\em Hardy-Littlewood $k$-tuple conjecture} \cite[p.\,62--68]{hl,riesel} 
suggests that admissible prime $k$-tuples near $x$ occur with average gaps
$a \approx C_k \log^k{\negthinspace}x$; the coefficients $C_k$ 
are reciprocal to the {\em Hardy-Littlewood $k$-tuple constants}\footnote{
For prime $k$-tuplets with $k\le7$, the numerical values are
$C_2\approx0.757392$,
$C_3\approx0.349864$,
$C_4\approx0.240895$,
$C_5\approx0.0986992$,  %4
$C_6\approx0.0578081$,  %1
$C_7\approx0.0185281$.  %4
%For prime $k$-tuplets with $k\ge8$, there may be more than one more constant $C_k$, depending on the specific $k$-tuplet pattern.
%For example, there are three different patterns of prime octuplets corresponding to two different values of $C_8$.
For definitions and further discussion, see also \cite{kourbatov2013}.
} 
\cite{finch,forbes}. 
Consider admissible prime $k$-tuples with a given pattern, and 
denote by $G_k(x)$ the largest gap between $k$-tuples below $x$. 
Then the following heuristic formula expresses 
the expected size of record gaps in terms of the average gap $a$, as $x\to\infty$:
$$
G_k(x) ~\sim~ a \log{x \over a} - ba ~ \lesssim ~  C_k \log^{k+1}{\negthinspace}x, 
 \quad \mbox{ with } \ a = C_k \log^k{\negthinspace}x, \ \ b \approx {2\over k}
 \qquad \mbox{\cite[p.\,5--7]{kourbatov2013}}.
$$
In the absence of rigorous upper bounds for gaps between prime constellations, 
it is important to accumulate extensive numerical evidence on actual record gaps.  
To that end, these tables contain all maximal gaps between prime {\em $k$-tuplets} (the densest admissible clusters of $k$ primes)
up to $10^{15}$, for each $k\le7$. 
Tables for $k=1$, 2, 4 have been previously published by several authors 
\cite{fischer,kourbatov2013,nicely,toes,rr,oeis}.
For $k=1$ the data extends up to $4\times10^{18}$ \cite{toes}.
Tables for $3\le k \le 7$ have been independently obtained as part of my own computation \cite{kourbatov2013,oeis}.
So is it {\bf\em always} true that $G_k(x)<C_k \log^{k+1}{\negthinspace}x$? 
It seems likely; no counterexamples thus far!

\newpage
\section{Record gaps between primes}

\begin{center}TABLE 1 \\Maximal gaps between primes below $4\times10^{18}$ \ 
\cite[OEIS A005250]{nicely,toes,oeis}
\\[0.5em]
\small
\begin{tabular}{rrr|rrr}  
\hline
\multicolumn{2}{r}{Consecutive primes}                            & Gap $g$ & 
\multicolumn{2}{r}{Consecutive primes~\phantom{\fbox{$11111^1$}}} & Gap $g$ \\
[0.5ex]\hline
\vphantom{\fbox{$1^1$}}
          2 &	3	    &   1 &        25056082087 &         25056082543 &  456 \\
          3 &	5	    &   2 &        42652618343 &         42652618807 &  464 \\
          7 &	11	    &   4 &       127976334671 &        127976335139 &  468 \\
         23 &	29	    &   6 &       182226896239 &        182226896713 &  474 \\
         89 &	97	    &   8 &       241160624143 &        241160624629 &  486 \\
        113 &	127	    &  14 &       297501075799 &        297501076289 &  490 \\
        523 &	541	    &  18 &       303371455241 &        303371455741 &  500 \\
        887 &	907	    &  20 &       304599508537 &        304599509051 &  514 \\
       1129 &	1151	    &  22 &       416608695821 &        416608696337 &  516 \\
       1327 &	1361	    &  34 &       461690510011 &        461690510543 &  532 \\
       9551 &	9587	    &  36 &       614487453523 &        614487454057 &  534 \\
      15683 &	15727	    &  44 &       738832927927 &        738832928467 &  540 \\
      19609 &	19661	    &  52 &      1346294310749 &       1346294311331 &  582 \\
      31397 &	31469	    &  72 &      1408695493609 &       1408695494197 &  588 \\
     155921 &	156007	    &  86 &      1968188556461 &       1968188557063 &  602 \\
     360653 &	360749	    &  96 &      2614941710599 &       2614941711251 &  652 \\
     370261 &	370373	    & 112 &      7177162611713 &       7177162612387 &  674 \\
     492113 &	492227	    & 114 &     13829048559701 &      13829048560417 &  716 \\
    1349533 &	1349651	    & 118 &     19581334192423 &      19581334193189 &  766 \\
    1357201 &	1357333	    & 132 &     42842283925351 &      42842283926129 &  778 \\
    2010733 &	2010881	    & 148 &     90874329411493 &      90874329412297 &  804 \\
    4652353 &	4652507	    & 154 &    171231342420521 &     171231342421327 &  806 \\
   17051707 &	17051887    & 180 &    218209405436543 &     218209405437449 &  906 \\
   20831323 &	20831533    & 210 &   1189459969825483 &    1189459969826399 &  916 \\
   47326693 &	47326913    & 220 &   1686994940955803 &    1686994940956727 &  924 \\
  122164747 &	122164969   & 222 &   1693182318746371 &    1693182318747503 & 1132 \\
  189695659 &	189695893   & 234 &  43841547845541059 &   43841547845542243 & 1184 \\
  191912783 &	191913031   & 248 &  55350776431903243 &   55350776431904441 & 1198 \\
  387096133 &	387096383   & 250 &  80873624627234849 &   80873624627236069 & 1220 \\
  436273009 &	436273291   & 282 & 203986478517455989 &  203986478517457213 & 1224 \\
 1294268491 &	1294268779  & 288 & 218034721194214273 &  218034721194215521 & 1248 \\
 1453168141 &	1453168433  & 292 & 305405826521087869 &  305405826521089141 & 1272 \\
 2300942549 &	2300942869  & 320 & 352521223451364323 &  352521223451365651 & 1328 \\
 3842610773 &	3842611109  & 336 & 401429925999153707 &  401429925999155063 & 1356 \\
 4302407359 &	4302407713  & 354 & 418032645936712127 &  418032645936713497 & 1370 \\
10726904659 &	10726905041 & 382 & 804212830686677669 &  804212830686679111 & 1442 \\
20678048297 &	20678048681 & 384 &1425172824437699411 & 1425172824437700887 & 1476 \\
22367084959 &	22367085353 & 394 &                    &                     &      \\
\hline
\end{tabular}
\end{center}
\normalsize

\newpage
\section{Record gaps between twin primes}

\begin{center}TABLE 2 \\Maximal gaps between twin primes $\{p$, $p+2\}$ \ 
\cite[OEIS A113274]{fischer,kourbatov2013,rr,oeis} \\[0.5em]
%\small
\begin{tabular}{rrr|rrr}  
\hline
\multicolumn{2}{r}{Initial primes in twin pairs}                        & Gap $g_2$ & 
\multicolumn{2}{r}{Initial primes in twin pairs~\phantom{\fbox{$1^1$}}} & Gap $g_2$ \\
[0.5ex]\hline
\vphantom{\fbox{$1^1$}}
               3 &               5 &             2 &     24857578817 &     24857585369 &          6552\\
               5 &              11 &             6 &     40253418059 &     40253424707 &          6648\\
              17 &              29 &            12 &     42441715487 &     42441722537 &          7050\\
              41 &              59 &            18 &     43725662621 &     43725670601 &          7980\\
              71 &             101 &            30 &     65095731749 &     65095739789 &          8040\\
             311 &             347 &            36 &    134037421667 &    134037430661 &          8994\\
             347 &             419 &            72 &    198311685749 &    198311695061 &          9312\\
             659 &             809 &           150 &    223093059731 &    223093069049 &          9318\\
            2381 &            2549 &           168 &    353503437239 &    353503447439 &         10200\\
            5879 &            6089 &           210 &    484797803249 &    484797813587 &         10338\\
           13397 &           13679 &           282 &    638432376191 &    638432386859 &         10668\\
           18539 &           18911 &           372 &    784468515221 &    784468525931 &         10710\\
           24419 &           24917 &           498 &    794623899269 &    794623910657 &         11388\\
           62297 &           62927 &           630 &   1246446371789 &   1246446383771 &         11982\\
          187907 &          188831 &           924 &   1344856591289 &   1344856603427 &         12138\\
          687521 &          688451 &           930 &   1496875686461 &   1496875698749 &         12288\\
          688451 &          689459 &          1008 &   2156652267611 &   2156652280241 &         12630\\
          850349 &          851801 &          1452 &   2435613754109 &   2435613767159 &         13050\\
         2868959 &         2870471 &          1512 &   4491437003327 &   4491437017589 &         14262\\
         4869911 &         4871441 &          1530 &  13104143169251 &  13104143183687 &         14436\\
         9923987 &         9925709 &          1722 &  14437327538267 &  14437327553219 &         14952\\
        14656517 &        14658419 &          1902 &  18306891187511 &  18306891202907 &         15396\\
        17382479 &        17384669 &          2190 &  18853633225211 &  18853633240931 &         15720\\
        30752231 &        30754487 &          2256 &  23275487664899 &  23275487681261 &         16362\\
        32822369 &        32825201 &          2832 &  23634280586867 &  23634280603289 &         16422\\
        96894041 &        96896909 &          2868 &  38533601831027 &  38533601847617 &         16590\\
       136283429 &       136286441 &          3012 &  43697538391391 &  43697538408287 &         16896\\
       234966929 &       234970031 &          3102 &  56484333976919 &  56484333994001 &         17082\\
       248641037 &       248644217 &          3180 &  74668675816277 &  74668675834661 &         18384\\
       255949949 &       255953429 &          3480 & 116741875898981 & 116741875918727 &         19746\\
       390817727 &       390821531 &          3804 & 136391104728629 & 136391104748621 &         19992\\
       698542487 &       698547257 &          4770 & 221346439666109 & 221346439686641 &         20532\\
      2466641069 &      2466646361 &          5292 & 353971046703347 & 353971046725277 &         21930\\
      4289385521 &      4289391551 &          6030 & 450811253543219 & 450811253565767 &         22548\\
     19181736269 &     19181742551 &          6282 & 742914612256169 & 742914612279527 &         23358\\
     24215097497 &     24215103971 &          6474 &1121784847637957 & 1121784847661339&         23382\\
\hline
\end{tabular}
\end{center}
\normalsize

\newpage
\section{Record gaps between prime triplets}

\begin{center}TABLE 3.1 \\Maximal gaps between prime triplets $\{p$, $p+2$, $p+6\}$ \ [OEIS A201598] \\[0.5em]
%\small
\begin{tabular}{rrr|rrr}  
\hline
\multicolumn{2}{r}{Initial primes in tuples}                        & Gap $g_3$ & 
\multicolumn{2}{r}{Initial primes in tuples~\phantom{\fbox{$1^1$}}} & Gap $g_3$ \\
[0.5ex]\hline
\vphantom{\fbox{$1^1$}}
        5 &        11 &      6 &      242419361 &       242454281 &  34920 \\
       17 &        41 &     24 &      913183487 &       913222307 &  38820 \\
       41 &       101 &     60 &     1139296721 &      1139336111 &  39390 \\
      107 &       191 &     84 &     2146630637 &      2146672391 &  41754 \\
      347 &       461 &    114 &     2188525331 &      2188568351 &  43020 \\
      461 &       641 &    180 &     3207540881 &      3207585191 &  44310 \\
      881 &      1091 &    210 &     3577586921 &      3577639421 &  52500 \\
     1607 &      1871 &    264 &     7274246711 &      7274318057 &  71346 \\
     2267 &      2657 &    390 &    33115389407 &     33115467521 &  78114 \\
     2687 &      3251 &    564 &    97128744521 &     97128825371 &  80850 \\
     6197 &      6827 &    630 &    99216417017 &     99216500057 &  83040 \\
     6827 &      7877 &   1050 &   103205810327 &    103205893751 &  83424 \\
    39227 &     40427 &   1200 &   133645751381 &    133645853711 & 102330 \\
    46181 &     47711 &   1530 &   373845384527 &    373845494147 & 109620 \\
    56891 &     58907 &   2016 &   412647825677 &    412647937127 & 111450 \\
    83267 &     86111 &   2844 &   413307596957 &    413307728921 & 131964 \\
   167621 &    171047 &   3426 &  1368748574441 &   1368748707197 & 132756 \\
   375251 &    379007 &   3756 &  1862944563707 &   1862944700711 & 137004 \\
   381527 &    385391 &   3864 &  2368150202501 &   2368150349687 & 147186 \\
   549161 &    553097 &   3936 &  2370801522107 &   2370801671081 & 148974 \\
   741677 &    745751 &   4074 &  3710432509181 &   3710432675231 & 166050 \\
   805031 &    809141 &   4110 &  5235737405807 &   5235737580317 & 174510 \\
   931571 &    937661 &   6090 &  8615518909601 &   8615519100521 & 190920 \\
  2095361 &   2103611 &   8250 & 10423696470287 &  10423696665227 & 194940 \\
  2428451 &   2437691 &   9240 & 10660256412977 &  10660256613551 & 200574 \\
  4769111 &   4778381 &   9270 & 11602981439237 &  11602981647011 & 207774 \\
  4938287 &   4948631 &  10344 & 21824373608561 &  21824373830087 & 221526 \\
 12300641 &  12311147 &  10506 & 36385356561077 &  36385356802337 & 241260 \\
 12652457 &  12663191 &  10734 & 81232357111331 &  81232357386611 & 275280 \\
 13430171 &  13441091 &  10920 &186584419495421 & 186584419772321 & 276900 \\
 14094797 &  14107727 &  12930 &297164678680151 & 297164678975621 & 295470 \\
 18074027 &  18089231 &  15204 &428204300934581 & 428204301233081 & 298500 \\
 29480651 &  29500841 &  20190 &450907041535541 & 450907041850547 & 315006 \\
107379731 & 107400017 &  20286 &464151342563471 & 464151342898121 & 334650 \\
138778301 & 138799517 &  21216 &484860391301771 & 484860391645037 & 343266 \\
156377861 & 156403607 &  25746 &666901733009921 & 666901733361947 & 352026 \\
\hline
\end{tabular}
%\medskip
%[OEIS A201598]
%
\end{center}
\normalsize

\newpage
\begin{center}TABLE 3.2  \\Maximal gaps between prime triplets $\{p$, $p+4$, $p+6\}$ \ [OEIS A201596]  % below $10^{15}$
\\[0.5em]
\small
\begin{tabular}{rrr|rrr}  
\hline
\multicolumn{2}{r}{Initial primes in tuples}                        & Gap $g_3$ & 
\multicolumn{2}{r}{Initial primes in tuples~\phantom{\fbox{$1^1$}}} & Gap $g_3$ \\
[0.5ex]\hline
\vphantom{\fbox{$1^1$}}
              7 &         13 &      6 &     2562574867 &      2562620653 &  45786 \\
             13 &         37 &     24 &     2985876133 &      2985923323 &  47190 \\
             37 &         67 &     30 &     4760009587 &      4760057833 &  48246 \\
            103 &        193 &     90 &     5557217797 &      5557277653 &  59856 \\
            307 &        457 &    150 &    10481744677 &     10481806897 &  62220 \\
            457 &        613 &    156 &    19587414277 &     19587476563 &  62286 \\
            613 &        823 &    210 &    25302582667 &     25302648457 &  65790 \\
           2137 &       2377 &    240 &    30944120407 &     30944191387 &  70980 \\
           2377 &       2683 &    306 &    37638900283 &     37638972667 &  72384 \\
           2797 &       3163 &    366 &    49356265723 &     49356340387 &  74664 \\
           3463 &       3847 &    384 &    49428907933 &     49428989167 &  81234 \\
           4783 &       5227 &    444 &    70192637737 &     70192720303 &  82566 \\
           5737 &       6547 &    810 &    74734558567 &     74734648657 &  90090 \\
           9433 &      10267 &    834 &   111228311647 &    111228407113 &  95466 \\
          14557 &      15643 &   1086 &   134100150127 &    134100250717 & 100590 \\
          24103 &      25303 &   1200 &   195126585733 &    195126688957 & 103224 \\
          45817 &      47143 &   1326 &   239527477753 &    239527584553 & 106800 \\
          52177 &      54493 &   2316 &   415890988417 &    415891106857 & 118440 \\
         126487 &     130363 &   3876 &   688823669533 &    688823797237 & 127704 \\
         317587 &     321817 &   4230 &   906056631937 &    906056767327 & 135390 \\
         580687 &     585037 &   4350 &   926175746857 &    926175884923 & 138066 \\
         715873 &     724117 &   8244 &  1157745737047 &   1157745878893 & 141846 \\
        2719663 &    2728543 &   8880 &  1208782895053 &   1208783041927 & 146874 \\
        6227563 &    6237013 &   9450 &  2124064384483 &   2124064533817 & 149334 \\
        8114857 &    8125543 &  10686 &  2543551885573 &   2543552039053 & 153480 \\
       10085623 &   10096573 &  10950 &  4321372168453 &   4321372359523 & 191070 \\
       10137493 &   10149277 &  11784 &  6136808604343 &   6136808803753 & 199410 \\
       18773137 &   18785953 &  12816 & 18292411110217 &  18292411310077 & 199860 \\
       21297553 &   21311107 &  13554 & 19057076066317 &  19057076286553 & 220236 \\
       25291363 &   25306867 &  15504 & 21794613251773 &  21794613477097 & 225324 \\
       43472497 &   43488073 &  15576 & 35806145634613 &  35806145873077 & 238464 \\
       52645423 &   52661677 &  16254 & 75359307977293 &  75359308223467 & 246174 \\
       69718147 &   69734653 &  16506 & 89903831167897 &  89903831419687 & 251790 \\
       80002627 &   80019223 &  16596 &125428917151957 & 125428917432697 & 280740 \\
       89776327 &   89795773 &  19446 &194629563521143 & 194629563808363 & 287220 \\
       90338953 &   90358897 &  19944 &367947033766573 & 367947034079923 & 313350 \\
      109060027 &  109081543 &  21516 &376957618687747 & 376957619020813 & 333066 \\
      148770907 &  148809247 &  38340 &483633763994653 & 483633764339287 & 344634 \\
     1060162843 & 1060202833 &  39990 &539785800105313 & 539785800491887 & 386574 \\
     1327914037 & 1327955593 &  41556 &  \multicolumn{3}{c}{}                     \\
\hline
\end{tabular}

%\medskip
%[OEIS A201596]
\end{center}
\normalsize

\newpage
\section{Record gaps between prime quadruplets}

\begin{center}TABLE 4 \\Maximal gaps between prime quadruplets $\{p$, $p+2$, $p+6$, $p+8\}$ \ [OEIS A113404]
\\[0.5em]
%\small
\begin{tabular}{rrr|rrr}  
\hline
\multicolumn{2}{r}{Initial primes in tuples}                        & Gap $g_4$ & 
\multicolumn{2}{r}{Initial primes in tuples~\phantom{\fbox{$1^1$}}} & Gap $g_4$ \\
[0.5ex]\hline
\vphantom{\fbox{$1^1$}}
         5 &         11 &      6 &      3043111031 &     3043668371 &  557340 \\
        11 &        101 &     90 &      3593321651 &     3593956781 &  635130 \\
       191 &        821 &    630 &      5675642501 &     5676488561 &  846060 \\
       821 &       1481 &    660 &     25346635661 &    25347516191 &  880530 \\
      2081 &       3251 &   1170 &     27329170151 &    27330084401 &  914250 \\
      3461 &       5651 &   2190 &     35643379901 &    35644302761 &  922860 \\
      5651 &       9431 &   3780 &     56390149631 &    56391153821 & 1004190 \\
     25301 &      31721 &   6420 &     60368686121 &    60369756611 & 1070490 \\
     34841 &      43781 &   8940 &     71335575131 &    71336662541 & 1087410 \\
     88811 &      97841 &   9030 &     76427973101 &    76429066451 & 1093350 \\
    122201 &     135461 &  13260 &     87995596391 &    87996794651 & 1198260 \\
    171161 &     187631 &  16470 &     96616771961 &    96618108401 & 1336440 \\
    301991 &     326141 &  24150 &    151023350501 &   151024686971 & 1336470 \\
    739391 &     768191 &  28800 &    164550390671 &   164551739111 & 1348440 \\
   1410971 &    1440581 &  29610 &    171577885181 &   171579255431 & 1370250 \\
   1468631 &    1508621 &  39990 &    210999769991 &   211001269931 & 1499940 \\
   2990831 &    3047411 &  56580 &    260522319641 &   260523870281 & 1550640 \\
   3741161 &    3798071 &  56910 &    342611795411 &   342614346161 & 2550750 \\
   5074871 &    5146481 &  71610 &   1970587668521 &  1970590230311 & 2561790 \\
   5527001 &    5610461 &  83460 &   4231588103921 &  4231591019861 & 2915940 \\
   8926451 &    9020981 &  94530 &   5314235268731 &  5314238192771 & 2924040 \\
  17186591 &   17301041 & 114450 &   7002440794001 &  7002443749661 & 2955660 \\
  21872441 &   22030271 & 157830 &   8547351574961 &  8547354997451 & 3422490 \\
  47615831 &   47774891 & 159060 &  15114108020021 & 15114111476741 & 3456720 \\
  66714671 &   66885851 & 171180 &  16837633318811 & 16837637203481 & 3884670 \\
  76384661 &   76562021 & 177360 &  30709975578251 & 30709979806601 & 4228350 \\
  87607361 &   87797861 & 190500 &  43785651890171 & 43785656428091 & 4537920 \\
 122033201 &  122231111 & 197910 &  47998980412211 & 47998985015621 & 4603410 \\
 132574061 &  132842111 & 268050 &  55341128536691 & 55341133421591 & 4884900 \\
 204335771 &  204651611 & 315840 &  92944027480721 & 92944033332041 & 5851320 \\
 628246181 &  628641701 & 395520 & 412724560672211 &412724567171921 & 6499710 \\
1749443741 & 1749878981 & 435240 & 473020890377921 &473020896922661 & 6544740 \\
2115383651 & 2115824561 & 440910 & 885441677887301 &885441684455891 & 6568590 \\
2128346411 & 2128859981 & 513570 & 947465687782631 &947465694532961 & 6750330 \\
2625166541 & 2625702551 & 536010 & 979876637827721 &979876644811451 & 6983730 \\
2932936421 & 2933475731 & 539310 &                 &                &         \\
\hline
\end{tabular}

%\medskip
%[OEIS A113404]
%
\end{center}
\normalsize

\newpage
\section{Record gaps between prime quintuplets}

\begin{center}TABLE 5.1 \\
Maximal gaps between prime quintuplets $\{p$, $p+2$, $p+6$, $p+8$, $p+12\}$  \\[0.5em]
%\small
\begin{tabular}{rrr|rrr}  
\hline
\multicolumn{2}{r}{Initial primes in tuples}                        & Gap $g_5$ & 
\multicolumn{2}{r}{Initial primes in tuples~\phantom{\fbox{$1^1$}}} & Gap $g_5$ \\
[0.5ex]\hline
\vphantom{\fbox{$1^1$}}
              5 &             11 &           6 &    107604759671 &   107616100511 &    11340840 \\
             11 &            101 &          90 &    140760439991 &   140772689501 &    12249510 \\
            101 &           1481 &        1380 &    162661360481 &   162673773671 &    12413190 \\
           1481 &          16061 &       14580 &    187735329491 &   187749510491 &    14181000 \\
          22271 &          43781 &       21510 &    327978626531 &   327994719461 &    16092930 \\
          55331 &         144161 &       88830 &    508259311991 &   508275672341 &    16360350 \\
         536441 &         633461 &       97020 &    620537349191 &   620554105931 &    16756740 \\
         661091 &         768191 &      107100 &    667672901711 &   667689883031 &    16981320 \\
        1461401 &        1573541 &      112140 &   1079628551621 &  1079646141851 &    17590230 \\
        1615841 &        1917731 &      301890 &   1104604933841 &  1104624218981 &    19285140 \\
        5527001 &        5928821 &      401820 &   1182148717481 &  1182168243071 &    19525590 \\
       11086841 &       11664551 &      577710 &   1197151034531 &  1197173264711 &    22230180 \\
       35240321 &       35930171 &      689850 &   2286697462781 &  2286720012251 &    22549470 \\
       53266391 &       54112601 &      846210 &   2435950632251 &  2435980618781 &    29986530 \\
       72610121 &       73467131 &      857010 &   3276773115431 &  3276805283951 &    32168520 \\
       92202821 &       93188981 &      986160 &   5229301162991 &  5229337555061 &    36392070 \\
      117458981 &      119114111 &     1655130 &   9196865051651 &  9196903746881 &    38695230 \\
      196091171 &      198126911 &     2035740 &  14660925945221 & 14660966101421 &    40156200 \\
      636118781 &      638385101 &     2266320 &  21006417451961 & 21006458070461 &    40618500 \\
      975348161 &      977815451 &     2467290 &  22175175736991 & 22175216733491 &    40996500 \\
     1156096301 &     1158711011 &     2614710 &  22726966063091 & 22727007515411 &    41452320 \\
     1277816921 &     1281122231 &     3305310 &  22931291089451 & 22931338667591 &    47578140 \\
     1347962381 &     1351492601 &     3530220 &  31060723328351 & 31060771959221 &    48630870 \\
     2195593481 &     2199473531 &     3880050 &  85489258071311 & 85489313115881 &    55044570 \\
     3128295551 &     3132180971 &     3885420 &  90913430825291 & 90913489290971 &    58465680 \\
     4015046591 &     4020337031 &     5290440 &  96730325054171 & 96730390102391 &    65048220 \\
     8280668651 &     8286382451 &     5713800 & 199672700175071 &199672765913051 &    65737980 \\
     9027127091 &     9033176981 &     6049890 & 275444947505591 &275445018294491 &    70788900 \\
    15686967971 &    15693096311 &     6128340 & 331992774272981 &331992848243801 &    73970820 \\
    18901038971 &    18908988011 &     7949040 & 465968834865971 &465968914851101 &    79985130 \\
    21785624291 &    21793595561 &     7971270 & 686535413263871 &686535495684161 &    82420290 \\
    30310287431 &    30321057581 &    10770150 & 761914822198961 &761914910291531 &    88092570 \\
\hline
\end{tabular}

\medskip
[OEIS A201073]
\end{center}
\normalsize

\newpage
\begin{center}TABLE 5.2 \\Maximal gaps between prime quintuplets $\{p$, $p+4$, $p+6$, $p+10$, $p+12\}$ \\[0.5em]
%\small
\begin{tabular}{rrr|rrr}  
\hline
\multicolumn{2}{r}{Initial primes in tuples}                        & Gap $g_5$ & 
\multicolumn{2}{r}{Initial primes in tuples~\phantom{\fbox{$1^1$}}} & Gap $g_5$ \\
[0.5ex]\hline
\vphantom{\fbox{$1^1$}}
              7 &             97 &          90 &     15434185927 &     15440743597 &     6557670\\
             97 &           1867 &        1770 &     17375054227 &     17381644867 &     6590640\\
           3457 &           5647 &        2190 &     17537596327 &     17544955777 &     7359450\\
           5647 &          15727 &       10080 &     25988605537 &     25997279377 &     8673840\\
          19417 &          43777 &       24360 &     66407160637 &     66416495137 &     9334500\\
          43777 &          79687 &       35910 &     74862035617 &     74871605947 &     9570330\\
         101107 &         257857 &      156750 &     77710388047 &     77723371717 &    12983670\\
        1621717 &        1830337 &      208620 &    144124106167 &    144138703987 &    14597820\\
        3690517 &        3995437 &      304920 &    210222262087 &    210238658797 &    16396710\\
        5425747 &        5732137 &      306390 &    585234882097 &    585252521167 &    17639070\\
        8799607 &        9127627 &      328020 &    926017532047 &    926036335117 &    18803070\\
        9511417 &        9933607 &      422190 &    986089952917 &    986113345747 &    23392830\\
       16388917 &       16915267 &      526350 &   2819808136417 &   2819832258697 &    24122280\\
       22678417 &       23317747 &      639330 &   3013422626107 &   3013449379477 &    26753370\\
       31875577 &       32582437 &      706860 &   3538026326827 &   3538053196957 &    26870130\\
       37162117 &       38028577 &      866460 &   4674635167747 &   4674662545867 &    27378120\\
       64210117 &       65240887 &     1030770 &   5757142722757 &   5757171559957 &    28837200\\
      119732017 &      120843637 &     1111620 &   7464931087717 &   7464961813867 &    30726150\\
      200271517 &      201418957 &     1147440 &   8402871269197 &   8402904566467 &    33297270\\
      203169007 &      204320107 &     1151100 &   9292699799017 &   9292733288557 &    33489540\\
      241307107 &      242754637 &     1447530 &  10985205390997 &  10985239010737 &    33619740\\
      342235627 &      344005297 &     1769670 &  12992848206847 &  12992884792957 &    36586110\\
      367358347 &      369151417 &     1793070 &  15589051692667 &  15589094176627 &    42483960\\
      378200227 &      380224837 &     2024610 &  24096376903597 &  24096421071127 &    44167530\\
      438140947 &      440461117 &     2320170 &  37371241083097 &  37371285854467 &    44771370\\
      446609407 &      448944487 &     2335080 &  38728669335607 &  38728728308527 &    58972920\\
      711616897 &      714020467 &     2403570 &  91572717670537 &  91572784840627 &    67170090\\
      966813007 &      970371037 &     3558030 & 109950817237357 & 109950886775827 &    69538470\\
     2044014607 &     2048210107 &     4195500 & 325554440818297 & 325554513360487 &    72542190\\
     3510456787 &     3514919917 &     4463130 & 481567288596127 & 481567361629087 &    73032960\\
     4700738167 &     4705340527 &     4602360 & 501796510663237 & 501796584764467 &    74101230\\
     5798359657 &     5803569847 &     5210190 & 535243109721577 & 535243185965557 &    76243980\\
     7896734467 &     7902065527 &     5331060 & 657351798174427 & 657351876771637 &    78597210\\
    12654304207 &    12659672737 &     5368530 & 818872754682547 & 818872840949077 &    86266530\\
    13890542377 &    13896088897 &     5546520 & 991851356676277 & 991851464273767 &   107597490\\
    14662830817 &    14668797037 &     5966220 &                 &                 &            \\
\hline
\end{tabular}

\medskip
[OEIS A201062]
\end{center}
\normalsize

\newpage
\section{Record gaps between prime sextuplets}
\begin{adjustwidth}{-9mm}{}
\begin{center}TABLE 6 \\Maximal gaps between prime sextuplets $\{p$, $p+4$, $p+6$, $p+10$, $p+12$, $p+16\}$ \\[0.5em]
\normalsize
\begin{tabular}{rrr|rrr}  
\hline
\multicolumn{2}{r}{Initial primes in tuples}                        & Gap $g_6$ & 
\multicolumn{2}{r}{Initial primes in tuples~\phantom{\fbox{$1^1$}}} & Gap $g_6$ \\
[0.5ex]\hline
\vphantom{\fbox{$1^1$}}
               7 &            97 &            90 &    422088931207 &   422248594837 &     159663630\\
              97 &         16057 &         15960 &    427190088877 &   427372467157 &     182378280\\
           19417 &         43777 &         24360 &    610418426197 &   610613084437 &     194658240\\
           43777 &       1091257 &       1047480 &    659829553837 &   660044815597 &     215261760\\
         3400207 &       6005887 &       2605680 &    660863670277 &   661094353807 &     230683530\\
        11664547 &      14520547 &       2856000 &    853633486957 &   853878823867 &     245336910\\
        37055647 &      40660717 &       3605070 &   1089611097007 &  1089869218717 &     258121710\\
        82984537 &      87423097 &       4438560 &   1247852774797 &  1248116512537 &     263737740\\
        89483827 &      94752727 &       5268900 &   1475007144967 &  1475318162947 &     311017980\\
        94752727 &     112710877 &      17958150 &   1914335271127 &  1914657823357 &     322552230\\
       381674467 &     403629757 &      21955290 &   1953892356667 &  1954234803877 &     342447210\\
      1569747997 &    1593658597 &      23910600 &   3428196061177 &  3428617938787 &     421877610\\
      2019957337 &    2057241997 &      37284660 &   9367921374937 &  9368397372277 &     475997340\\
      5892947647 &    5933145847 &      40198200 &  10254799647007 & 10255307592697 &     507945690\\
      6797589427 &    6860027887 &      62438460 &  13786576306957 & 13787085608827 &     509301870\\
     14048370097 &   14112464617 &      64094520 &  21016714812547 & 21017344353277 &     629540730\\
     23438578897 &   23504713147 &      66134250 &  33157788914347 & 33158448531067 &     659616720\\
     24649559647 &   24720149677 &      70590030 &  41348577354307 & 41349374379487 &     797025180\\
     29637700987 &   29715350377 &      77649390 &  72702520226377 & 72703333384387 &     813158010\\
     29869155847 &   29952516817 &      83360970 &  89165783669857 & 89166606828697 &     823158840\\
     45555183127 &   45645253597 &      90070470 & 122421000846367 &122421855415957 &     854569590\\
     52993564567 &   53086708387 &      93143820 & 139864197232927 &139865086163977 &     888931050\\
     58430706067 &   58528934197 &      98228130 & 147693859139077 &147694869231727 &    1010092650\\
     93378527647 &   93495691687 &     117164040 & 186009633998047 &186010652137897 &    1018139850\\
     97236244657 &   97367556817 &     131312160 & 202607131405027 &202608270995227 &    1139590200\\
    240065351077 &  240216429907 &     151078830 & 332396845335547 &332397997564807 &    1152229260\\
    413974098817 &  414129003637 &     154904820 & 424681656944257 &424682861904937 &    1204960680\\
    419322931117 &  419481585697 &     158654580 & 437804272277497 &437805730243237 &    1457965740\\
\hline
\end{tabular}

\medskip
[OEIS A200503]
\end{center}
\end{adjustwidth}
\normalsize

\phantom{\vspace{1cm}}

\newpage
\section{Record gaps between prime septuplets}

\begin{adjustwidth}{-9mm}{}
\begin{center}TABLE 7.1 
\\Maximal gaps between prime 7-tuples $\{p$, $p+2$, $p+8$, $p+12$, $p+14$, $p+18$, $p+20\}$ \\[0.5em]
%\small
\begin{tabular}{rrr|rrr}  
\hline
\multicolumn{2}{r}{Initial primes in tuples}                        & Gap $g_7$ & 
\multicolumn{2}{r}{Initial primes in tuples~\phantom{\fbox{$1^1$}}} & Gap $g_7$ \\
[0.5ex]\hline
\vphantom{\fbox{$1^1$}}
           5639 &          88799 &       83160 &   1554893017199 &  1556874482069 &  1981464870\\
          88799 &         284729 &      195930 &   2088869793539 &  2090982626639 &  2112833100\\
         284729 &         626609 &      341880 &   2104286376329 &  2106411289049 &  2124912720\\
        1146779 &        6560999 &     5414220 &   2704298257469 &  2706823007879 &  2524750410\\
        8573429 &       17843459 &     9270030 &   3550904257709 &  3553467600029 &  2563342320\\
       24001709 &       42981929 &    18980220 &   4438966968419 &  4442670730019 &  3703761600\\
       43534019 &       69156539 &    25622520 &   9996858589169 & 10000866474869 &  4007885700\\
       87988709 &      124066079 &    36077370 &  21937527068909 & 21942038052029 &  4510983120\\
      157131419 &      208729049 &    51597630 &  29984058230039 & 29988742571309 &  4684341270\\
      522911099 &      615095849 &    92184750 &  30136375346249 & 30141383681399 &  5008335150\\
      706620359 &      832143449 &   125523090 &  32779504324739 & 32784963061379 &  5458736640\\
     1590008669 &     1730416139 &   140407470 &  40372176609629 & 40377635870639 &  5459261010\\
     2346221399 &     2488117769 &   141896370 &  42762127106969 & 42767665407989 &  5538301020\\
     3357195209 &     3693221669 &   336026460 &  54620176867169 & 54626029928999 &  5853061830\\
    11768282159 &    12171651629 &   403369470 &  63358011407219 & 63365153990639 &  7142583420\\
    30717348029 &    31152738299 &   435390270 &  79763188368959 & 79770583970249 &  7395601290\\
    33788417009 &    34230869579 &   442452570 & 109974651670769 &109982176374599 &  7524703830\\
    62923039169 &    63550891499 &   627852330 & 145568747217989 &145576919193689 &  8171975700\\
    68673910169 &    69428293379 &   754383210 & 196317277557209 &196325706400709 &  8428843500\\
    88850237459 &    89858819579 &  1008582120 & 221953318490999 &221961886287509 &  8567796510\\
   163288980299 &   164310445289 &  1021464990 & 249376874266769 &249385995968099 &  9121701330\\
   196782371699 &   197856064319 &  1073692620 & 290608782523049 &290618408585369 &  9626062320\\
   421204876439 &   422293025249 &  1088148810 & 310213774327979 &310225023265889 & 11248937910\\
   427478111309 &   428623448159 &  1145336850 & 471088826892779 &471100312066829 & 11485174050\\
   487635377219 &   489203880029 &  1568502810 & 631565753063879 &631578724265759 & 12971201880\\
   994838839439 &   996670266659 &  1831427220 & 665514714418439 &665530090367279 & 15375948840\\
\hline
\end{tabular}

\medskip
[OEIS A201251]
\end{center}
\end{adjustwidth}
\normalsize

\newpage
\begin{adjustwidth}{-9mm}{}
\begin{center}TABLE 7.2 \\
Maximal gaps between prime 7-tuples $\{p$, $p+2$, $p+6$, $p+8$, $p+12$, $p+18$, $p+20\}$  \\[0.5em]
%\small
\begin{tabular}{rrr|rrr}  
\hline
\multicolumn{2}{r}{Initial primes in tuples}                        & Gap $g_7$ & 
\multicolumn{2}{r}{Initial primes in tuples~\phantom{\fbox{$1^1$}}} & Gap $g_7$ \\
[0.5ex]\hline
\vphantom{\fbox{$1^1$}}
             11 &         165701 &     165690 &   382631592641 &   383960791211 & 1329198570\\
         165701 &        1068701 &     903000 &   711854781551 &   714031248641 & 2176467090\\
        1068701 &       11900501 &   10831800 &  2879574595811 &  2881987944371 & 2413348560\\
       25658441 &       39431921 &   13773480 &  3379186846151 &  3381911721101 & 2724874950\\
       45002591 &       67816361 &   22813770 &  5102247756491 &  5105053487531 & 2805731040\\
       93625991 &      124716071 &   31090080 &  5987254671311 &  5990491102691 & 3236431380\\
      257016491 &      300768311 &   43751820 &  7853481899561 &  7857040317011 & 3558417450\\
      367438061 &      428319371 &   60881310 & 11824063534091 & 11828142800471 & 4079266380\\
      575226131 &      661972301 &   86746170 & 16348094430581 & 16353374758991 & 5280328410\\
     1228244651 &     1346761511 &  118516860 & 44226969237161 & 44233058406611 & 6089169450\\
     1459270271 &     1699221521 &  239951250 & 54763336591961 & 54771443197181 & 8106605220\\
     2923666841 &     3205239881 &  281573040 &154325181803321 &154333374270191 & 8192466870\\
    10180589591 &    10540522241 &  359932650 &157436722520921 &157445120715341 & 8398194420\\
    15821203241 &    16206106991 &  384903750 &281057032201481 &281065611322031 & 8579120550\\
    23393094071 &    23911479071 &  518385000 &294887168565161 &294896169845351 & 9001280190\\
    37846533071 &    38749334621 &  902801550 &309902902299701 &309914040972071 &11138672370\\
   158303571521 &   159330579041 & 1027007520 &419341934631071 &419354153220461 &12218589390\\
   350060308511 &   351146640191 & 1086331680 &854077393259801 &854090557643621 &13164383820\\
\hline
\end{tabular}

\medskip
[OEIS A201051]
\end{center}
\end{adjustwidth}

\bigskip
\hrule
\bigskip

\noindent 2010 {\it Mathematics Subject Classification}: 11N05
%Primary 11N05; Secondary 60G70.

\noindent \emph{Keywords: } 
distribution of primes, prime $k$-tuple, Hardy-Littlewood conjecture,
extreme value statistics, prime gap, Cram\'er conjecture, 
prime constellation, twin prime conjecture, prime triplet, prime quadruplet, 
prime quintuplet, prime sextuplet, prime septuplet.

\bigskip
\hrule
\bigskip

\noindent (Concerned with OEIS sequences
A005250,
A113274,
A113404,
A200503,
A201051,
A201062,
A201073,
A201251,
A201596,
A201598,
A202281,
A202361.)

\bigskip

\end{document}